\documentclass[12pt]{article}

\usepackage{amssymb}
\usepackage{fullpage}

\newcommand{\col}{\mathop{\mathrm{col}}\limits}

\newcommand{\ON}{\mathrm{ON}}

\newtheorem{definition}{Definition}
\newcommand{\Ubf}{\widetilde}
\newtheorem{lemma}[definition]{Lemma}
\newtheorem{claim}[definition]{Claim}
\newcommand{\proof}{\noindent{\bf Proof.}\ \ }
\newcommand{\proofend}{\hfill~$\square$}
\newcommand{\crit}{\mathop{\mathrm{crit}}\nolimits}
\newtheorem{theorem}[definition]{Theorem}
\newtheorem{quotethm}{Theorem}
\newtheorem{quotelemma}{Lemma}
\newtheorem{quoteclaim}{Claim}

\let\<=\langle
\let\>=\rangle
\let\cross=\times
\newcommand{\into}{\rightarrow}
\newcommand{\so}{\{\,}
\newcommand{\sr}{\,\}}
\newcommand{\ult}[2]{\mathop{\mathrm{Ult}}\limits(#1,#2)}

\newcommand{\rest}{\restriction} 	
\newcommand{\R}{\mathord{\mathbb R}}

\newcommand{\forces}{\Vdash}

\title{Proper forcing and $L(\R)$} 

\author{Itay Neeman \thanks{Partially supported by NSF grant DMS 98-03292.}
\\ Harvard University \\ Department of Mathematics \\
Cambridge, MA 02138 \\ ineeman@math.harvard.edu
\and Jind\v rich Zapletal \thanks{The second author 
acknowledges support from NSF grant DMS 9022140, GA \v CR grant 201/97/0216
and CRM, Universita Aut\' onoma de Barcelona.}\\ Department of Mathematics \\
Dartmouth College \\ Hanover, NH 03755 \\ Jindrich.Zapletal@dartmouth.edu}

\begin{document}

\maketitle

\begin{abstract}
We present two ways in which the model $L(\R)$ is canonical assuming the existence of 
large cardinals. We show that the theory of this model, with {\em ordinal} 
parameters, cannot be changed by small forcing; we show further that a set of 
ordinals in $V$ cannot be added to $L(\R)$ by small forcing. The large cardinal 
needed corresponds to the consistency strength of $AD^{L(\R)}$; roughly $\omega$
Woodin cardinals.
\end{abstract}

\section{Introduction}

It is well known that under the existence of large cardinals
the theory of $L(\R)$ ---possibly with real parameters--- is absolute, and in 
particular cannot be changed by small forcings. Things may be different if one 
considers the theory of $L(\R)$ with {\em ordinal} parameters. 
By results
of Woodin and Shelah this theory can be changed by semiproper 
forcing, even granted large cardinals. In fact the truth value of the formula
$\phi(\Ubf{\alpha})=\mbox{``$\Ubf{\alpha}$ is equal to $\omega_2$''}$ 
can easily be changed as follows:
Start with $V$ having, e.g., a supercompact.
Then $L(\R^V)\models \phi[\alpha]$ for some 
$\alpha\leq\omega_2^V$, simply because any cardinal of $V$ must
be a cardinal of $L(\R^V)$.
Using the supercompact force to make the SemiProper Forcing Axiom hold in 
the generic extension. The supercompact of $V$ becomes $\omega_2$ of $V[G]$, and 
in $V[G]$ every set still has a sharp. All this can be done with a semiproper 
forcing.
From SPFA it follows in $V[G]$ that the non-stationary 
ideal on $\omega_1$ is 
saturated (see \cite{Jech-MF}; this result is due to Shelah), 
and by results of Woodin this (with sharps) implies
that $\omega_2$ as computed in $L(\R^{V[G]})$ is equal
to $\omega_2^{V[G]}$ (see \cite{Woodin-NS-ideal}). Now this is greater than
$\omega_2^V$, and so certainly greater than $\alpha$. 
Thus $L(\R^V)\models\phi[\alpha]$ while $L(\R^{V[G]})\not\models\phi[\alpha]$.

This example demonstrates that semiproper forcings {\em can} change the theory of 
$L(\R)$ with ordinal parameters greater than or equal to $\omega_2$.\footnote{Ordinal 
parameters below $\omega_2$ can, modulo $\omega_1$, be coded by reals. Thus, assuming 
large cardinals, forcing notions which preserve $\omega_1$ cannot change the theory 
of $L(\R)$ with ordinal parameters below $\omega_2$.} Any attempt to prove the 
preservation of this theory must therefore be restricted to a class of forcing 
strictly smaller than semiproper. 

\begin{theorem}[Embedding Theorem]
(Under large cardinal assumption $A_\kappa$, see below.)
Let $P$ be a proper forcing notion of size $\leq\kappa$, and let $G$ be
$P$-generic/$V$. Then there exists an elementary embedding 
$j\colon L(\R^V) \into L(\R^{V[G]})$ which is the identity on all ordinals.
\label{embedding}
\end{theorem}
The large cardinal assumption, $A_\kappa$, of Theorem \ref{embedding} is the 
following:

\vspace{0.1in}
\noindent
$(A_\kappa)$ \hspace{0.125in}
There exists a class inner model $M$ and a countable ordinal $\delta$ so
that
\begin{itemize}
\item $M=L(V_\delta^M)$;
\item $M\models\mbox{``$\delta$ is the supremum of $\omega$ Woodin 
cardinals;''}$ and
\item $M$ is uniquely iterable for iteration trees of length $\leq\kappa^+$.
\end{itemize}
{\em Uniquely iterable} means basically that in the course of our proof
we are free to create iteration trees on $M$, without having to worry about
the existence of cofinal well-founded branches.
More precisely $M$ must be iterable, meaning that the good player must win the 
full iteration game of \cite{it} on $M$ of length $\kappa^++1$. Furthermore the 
choice of cofinal branches must be unique, in the sense that for any iteration 
tree on $M$ of size $<\kappa^+$ there must be a {\em unique} cofinal branch $b$ 
such that the direct limit model $M_b$ is itself iterable. The technical 
assumption $A_\kappa$ is weaker than the existence (above $\kappa$) of $\omega$ 
Woodin cardinals and a measurable cardinal above them. It is closely connected 
to the large cardinal strength of $AD^{L(\R)}$.

Theorem \ref{embedding} implies in particular that the example given above 
cannot be carried out with proper (as opposed to semiproper) forcing; 
the full theory of $L(\R)$, with ordinal parameters, cannot be changed by proper 
forcing. It is a further immediate corollary of the Embedding Theorem
that ${\mathrm {HOD}}^{L(\R)}$ cannot be changed by proper forcing. Proper 
forcings also cannot ``code'' into $L(\R)$ a set of ordinals $A\in 
V\setminus L(\R)$: 

\begin{theorem}[Anti-coding Theorem]
(Under large cardinal assumption $A_\kappa$). Let $P$ be a proper
forcing notion of size $\leq\kappa$, and let $G$ be $P$-generic/$V$.
Suppose that $A\subset\ON$ is in $V$; 
then $A\in L(\R^V)\iff  A\in L(\R^{V[G]})$.
\label{anti-coding}
\end{theorem}

As with the Embedding Theorem, the Anti-Coding Theorem
cannot be extended much further. By a result of Woodin it fails
for semiproper forcings (provably from large cardinals). 
Both the Anti-Coding Theorem and the Embedding Theorem
do however extend to the class of reasonable forcings --- a class slightly 
bigger than proper. The proofs in this paper apply to reasonable forcings.

Our Theorems are similar in flavor to 
results of Foreman and Magidor
\cite{FM}, who investigated the possibility of forcing to change the
{\em definable continuum} --- the supremum of all ordinals $\gamma$
such that $\gamma$ is the order type of some prewellorderings
of reals in $L(\R)$. (This ordinal is commonly denoted as $\theta^{
L(\R)}$).
In \cite{FM} it is shown that (granted large cardinals) reasonable
forcings cannot change the definable continuum. This result can be obtained
also from our Embedding Theorem, since clearly
$\theta^{L(\R^V)}=j(\theta^{L(\R^V)})=\theta^{L(\R^{V[G]})}$
where $G$ is $P$-generic/$V$ for a reasonable forcing $P$, and
$j$ is the elementary embedding given by Theorem \ref{embedding}.
\cite{FM} prove a more general result concerning prewellorderings which
are homogeneously Suslin. Using additional work of Woodin's it is 
possible to derive the Embedding Theorem from their result. The proof
of the Embedding Theorem which we include here is different, and
its methods are needed later to obtain the Anti-Coding Theorem.

\vspace{0.1in}

As with many results involving $L(\R)$ and large cardinals there are
(at least)
two alternative routes to proving our Theorems; one which uses stationary
tower forcing, and another which uses iteration trees. The latter is presented 
in this paper while the former can be found in \cite{NZ-stf}. It is interesting 
that even though iteration trees and stationary tower forcing
are technically entirely different there are several similarities 
between the two approaches. Historically the Embedding and Anti-Coding Theorems 
were conjectured by the second author, who from a weakly compact
Woodin cardinal proved the first for c.c.c forcings and
the second for c.c.c. forcing as well as proper forcing
notions contained in $\omega_1$. Both proofs used the techniques
of stationary tower forcing. Those results were presented
during the 1996 Set Theory meeting in Luminy, France. The first
author subsequently used iteration trees to prove the full
Theorems as they appear in the present paper, while the
second author strengthened the stationary tower proofs to prove roughly the same
results as they appear in \cite{NZ-stf}.

The structure of this paper is such that most of the use of large cardinal
assumptions is exiled into two ``black boxes'' (Woodin's genericity
iterations) which are quoted and then used. The proofs relating to these black 
boxes are due entirely to Hugh Woodin. Readers who are not
experts on large cardinals may still be able to follow the proofs of the 
Embedding and Anti-Coding Theorems if they are willing to accept these black 
boxes. In Section \ref{sec.embedding}
we present the proof of the Embedding Theorem, and in Section
\ref{sec.anti-coding} the proof of the Anti-Coding Theorem. The proof in Section 
\ref{sec.anti-coding} uses the techniques of Section \ref{sec.embedding} as its 
backbone. The proofs of the black boxes are included in an Appendix to the
e-print of this paper at http://arXiv.org.

\section{The Embedding Theorem}

\label{sec.embedding}
 
We begin now the proof of the Embedding Theorem. Fix a proper
forcing notion $P$ and a generic $G$. Fix $M$ which witnesses $A_\kappa$. 
To prove the Theorem we must
construct the elementary embedding $j\colon L(\R^V) \into L(\R^{V[G]})$.
The requirement that $j\rest \ON$ be the identity essentially tells us
what $j$ is. We must have $j(x)=x$ for any $x\in\R^V$, and since 
all elements of $L(\R)$ are definable from a real and some ordinals
this fixes the map $j$ completely. Any element of $L(\R^V)$ definable
in $L(\R^V)$ from the real $z$ and the ordinals $\alpha_0,\dots\alpha_k$
using the formula $\phi$, 
must be mapped to the element of $L(\R^{V[G]})$ definable from 
$z,\alpha_0,\dots,\alpha_k$ using the same formula $\phi$ in
$L(\R^{V[G]})$. All we must prove is that this gives a well-defined
elementary embedding $j$, and this amounts to showing that
\begin{center}
\parbox{4.2in}{
For any $z\in\R^V$, any $\alpha_0,\dots,\alpha_k\in\ON$,
and any formula $\phi$, 
$$
\begin{array}{l}
L(\R^V)\models\phi[z,\alpha_0,\dots,\alpha_k] \iff \\
	L(\R^{V[G]})\models\phi[z,\alpha_0,\dots,\alpha_k].
\end{array}$$
}
\end{center}
Fix $z$, $\alpha_0,\dots,\alpha_k$, and a formula 
$\phi$. We shall prove the above equivalence using a symmetric collapse. Given a 
model $N$ 
and some ordinal $\lambda$, we consider the L\'evy Collapse
$\col(\omega,<\!\lambda)$ ---
the finite support product of the forcings $\col(\omega,\xi)$ for 
$\xi<\lambda$. Define the name
$\dot{\R}_{symm}=\bigcup_{\xi<\lambda}\R^{N[\dot{H}\rest \col(\omega,<\xi)]}$,
where $\dot{H}$ is a name for the generic object. 
$\dot{\R}_{symm}$ are the reals in the symmetric collapse up to 
$\lambda$. Those were first investigated by Solovay who used a symmetric
collapse to construct a model where all sets of reals are Lebesgue measurable.
The important property of the collapse is its
{\em homogeneity} --- any statement about $\dot{\R}_{symm}$ which involves only 
parameters from $N$ is true in the generic extension iff it is forced by the 
{\em empty condition} (see \cite{Jech}).

Our strategy is to construct a model $N$ and two different
generics $H_1$ and $H_2$ such that
\begin{enumerate}
\item $z\in N$;
\item $H_1$ and $H_2$ are both $\col(\omega,<\!\lambda)$-generic/$N$;
\item $\dot{\R}_{symm}[H_1]=\R^V$; and 
\item $\dot{\R}_{symm}[H_2]=\R^{V[G]}$.
\end{enumerate}
This will immediately complete the proof, as
\begin{eqnarray*}
L(\R^V)\models\phi[z,\alpha_0,\dots,\alpha_k]	& \iff_1 & 
N[H_1]\models\mbox{``$L(\dot{\R}_{symm}[H_1])\models
		\phi[z,\alpha_0,\dots,\alpha_k]$''} \\
	& \iff_2 & N\models\mbox{``$\forces^{\col(\omega,<\!\lambda)}
L(\dot{\R}_{symm})\models\phi[z,\alpha_0,\dots,\alpha_k]$''} \\
	& \iff_3 & N[H_2]\models\mbox{``$L(\dot{\R}_{symm}
[H_2])\models\phi[z,\alpha_0,\dots,\alpha_k]$''} \\
 	& \iff_4 & L(\R^{V[G]})\models\phi[z,\alpha_0,\dots,\alpha_k].
\end{eqnarray*}
The implications $1$ and $4$ follow from items
(3) and (4) above. The implications $2$ and $3$ follow from the homogeneity 
of the forcing.

We construct $N$ as an iterate of the model $M$, in $\omega$ stages. 
Each stage will be carried out in $V$ while the full construction
will exist in $V^{\col(\omega,\R)}$. Our main tool
is the following Theorem of Woodin's (see \cite{stf} or 
http://www.???.???).

\begin{quotethm}[Woodin's first genericity iteration]
Let $Q$ be an $\omega_1+1$-iterable inner model and 
let $\tau<\eta$ be countable 
(in $V$) ordinals such
that $Q\models\mbox{``$\eta$ is a Woodin cardinal''}$.
{\bf Then} there exists a forcing notion ${\mathbb W}^Q_{\tau,\eta}\in Q$ 
of size $\eta$,
such that for {\em any} real $x$ it is possible to construct an
iteration embedding $j\colon Q \into \tilde{Q}$
with the property that 
\begin{itemize}
\item $x$ is $j({\mathbb W}^Q_{\tau,\eta})$-generic/$\tilde{Q}$;
\item $j(\eta)$ is countable in $V$, indeed $j''(\omega_1^V)\subset
							\omega_1^V$; and
\item $\crit(j)>\tau$.
\end{itemize}
Furthermore for any small forcing ${\mathbb O}\in V_\tau^Q$
there exists an $\mathbb O$ name for a forcing notion 
$\dot{\mathbb W}^{Q,\mathbb O}_{\tau,\eta}$ so that
for any $o$ which is
$\mathbb O$-generic/$Q$, there exists an iteration 
embedding $j\colon Q\into \tilde{Q}$ satisfying the above except that
now $x$ is made 
$j(\dot{\mathbb W}^{Q,\mathbb O}_{\tau,\eta})[o]$-generic/$\tilde{Q}[o]$.\footnote{
This does not follow from the previous part of the Theorem,
since it gives a $j$ which is an iteration of $Q$, and this is 
more restrictive than being an iteration of $Q[o]$.}
\end{quotethm}

Woodin's genericity Theorem immediately tells us how to iterate
$M$ so as to satisfy condition (3) above. Fix some
$g\colon \omega \into \R^V$ which is $\col(\omega,\R)$-generic/$V$,
and so enumerates all the reals of $V$.
Our plan is to apply Woodin's Theorem using
the $2i$-th Woodin cardinal of $M$ to make $g(i)$ generic
over an iterate of $M$. (The reason we use only the even Woodin cardinals
will become clear presently.)
More precisely, working in $M$ we let $\so\delta_i\sr_{i\in\omega}$
be an increasing sequence of Woodin cardinals with supremum $\delta$.
Inductively define $\dot{\mathbb B}_i$ to be Woodin's forcing 
$\dot{\mathbb W}^{M,\dot{\mathbb B}_0\ast\dots\ast\dot{\mathbb B}_{i-1}}_{
\delta_{2i-1},\delta_{2i}}$ (defined in $M$).
Let $\mathbb B$ be the finite support iteration of the forcings
$\dot{\mathbb B}_i$.

Inductively construct an iteration of $M$. Begin by letting $M_0=M$, and
construct models $M_i$ and embeddings $j_i\colon M_i \into M_{i+1}$ so that
\begin{itemize}
\item[a.] $\< g(0),\dots,g(i-1)\>$ is $j_{0,2i}(\dot{\mathbb B}_0\ast\dots\ast
	\dot{\mathbb B}_{i-1})$-generic over $M_{2i}$, where $j_{0,2i}$ is 
	obtained
	through composition of the $j_i$'s.
\item[b.] $j_{2i}\colon M_{2i} \into M_{2i+1}$ is an embedding 
	to make the real $g(i)$ generic
	for the forcing notion
	$j_{0,2i+1}(\dot{\mathbb B}_i)[g(0),\dots,g(i-1)]$
	over the model $M_{2i+1}[g(0),\dots,g(i-1)]$. (Such an embedding
	always exists by Woodin's first genericity iteration.)
	We take $j_{2i}$ to be the identity whenever possible.
\item[c.] For the time being, let
	$j_{2i+1}\colon M_{2i+1}\into M_{2i+2}$ be the identity.
\end{itemize}

When iterating $M$ we use the unique iteration strategy. Thus by
an ``iteration embedding of $M$'' we mean only embeddings obtained
through those iteration trees on $M$ which choose the unique
iterable branch at every limit stage. By our iterability assumption
on $M$ this guarantees that direct limits of iteration 
embeddings in $V$ are well founded. Let $M_\infty$ be the direct limit model of the 
$M_i$-s and let $j_{i,\infty}\colon M_i \into M_\infty$
be the direct limit maps. Observe that $M_\infty$ is well founded. This is not 
entirely trivial as the {\em sequence} $j_i$ does not belong to $V$. However if this 
sequence gave rise to an ill founded direct limit one could use Schoenfield 
absoluteness to pull the existence of such a ``bad'' sequence back to $V$. Notice 
further that $j_{2i,\infty}$ has critical point greater 
than
$j_{0,2i}(\delta_{2i-1})$ so that
$j_{0,\infty}( \dot{\mathbb B}_0\ast\dots\ast
	 \dot{\mathbb B}_{i-1})=j_{0,2i}( \dot{\mathbb B}_0\ast\dots\ast
	 \dot{\mathbb B}_{i-1})$ and 
$j_{0,\infty}(\delta_{2i-1})=j_{0,{2i}}(\delta_{2i-1})$. In particular, 
$j_{0,\infty}(\delta_{2i-1})$ is countable (in $V$) and so
$j_{0,\infty}(\delta)\leq\omega_1^V$. 

It is clear from the construction that $\<g(0),\dots,g(i-1)\>$ is
$j_{0,\infty}( \dot{\mathbb B}_0\ast\dots\ast
	 \dot{\mathbb B}_{i-1})$-generic/$M_\infty$ for all $i$. Using the fact 
that $g$ is $\col(\omega,\R^V)$-generic/$V$ one can verify further that 
$\<g(i)\:\mid\: i<\omega\>$ is 
$j_{0,\infty}({\mathbb B})$-generic over $M_\infty$. Specifically,
fix any dense set $D$ in $j_{0,\infty}({\mathbb B})$ and
assume $(*)$ that the filter given by $\<g(i)\:\mid\: i<\omega\>$ does not
intersect $D$. Fix some $n$ large enough so that
$\<g(0),\dots,g(n-1)\>$ forces $(*)$ in $\col(\omega,\R)$, and such that
$D=j_{2n,\infty}(\bar{D})$ for some
$\bar{D}\in M_{2n}$. Now by condition (a), $\<g(0),\dots,g(n-1)\>$
is $j_{0,2n}(\dot{\mathbb B}_0\ast\dots\ast\dot{
\mathbb B}_{n-1})$-generic/$M_{2n}$. Working in $M_{2n}[g(0),\dots,g(n-1)]$ we 
can therefore find a condition 
$b=\<\dot{b}_0,\dots,\dot{b}_{k-1}\>\in \bar{D}$ 
(with $k\geq n$) such that $\<\dot{b}_0,\dots,\dot{b}_{n-1}\>$ belongs to the
$j_{0,2n}({\mathbb B}_0\ast\dots\ast{\mathbb B}_{n-1})$-generic given
by $\<g(0),\dots,g(n-1)\>$.
Next let us force over $M_{2n}$ with $j_{0,2n}(\dot{\mathbb B}_n\ast\dots\ast
\dot{\mathbb B}_{k-1})[g(0),
\dots,g(n-1)]$ below the condition $\<\dot{b}_n,\dots,\dot{b}_{k-1}\>[
g(0),\dots,g(n-1)]$,
and obtain reals $y_n,\dots,y_{k-1}$ such that 
\begin{itemize}
\item[G1.] $\<g(0),\dots,g(n-1),y_n,\dots,y_{k-1}\>$ is $j_{0,2n}(\dot{\mathbb 
B}_0\ast
\dots\ast\dot{\mathbb B}_{k-1})$-generic/$M_{2n}$, and
\item[G2.] this generic contains the condition $b\in\bar{D}$.
\end{itemize}
Such reals can be found in $V$ since the level of $M_{2n}$ involved is countable 
in $V$ (see conditions (i,ii) below). Consider finally the condition 
$\<g(0),\dots,g(n-1),
y_n,\dots,y_{k-1}\>$ in the forcing $\col(\omega,\R)$. This condition forces  
our construction to produce a model $M_{2k}$ which is equal to $M_{2n}$, and an 
embedding $j_{2n,2k}$ which equals the identity (note our use here of the requirement 
in (b) that $j_{2i}$ be the identity whenever possible). 
From this together with (G1,G2) 
it follows easily that $(*)$ is forced to {\em fail}, but this is a contradiction.

Observe next that the forcing $\mathbb B$ can be replaced by a symmetric
collapse. 
In other words it is possible to find 
$H\in M_\infty[g(i)\:\mid\: i<\omega]$ which is 
$\col(\omega,<\!j_{0,\infty}(\delta))$-generic/$M_\infty$ and so 
that $\dot{\R}_{symm}^{M_\infty}[H]=
\so g(i)\:\mid\: i<\omega \sr$. 
In fact it is well known that in general 
(for $\delta$ a strong limit cardinal) whenever
$\mathbb A$ is a direct limit of a regular 
chain of forcings ${\mathbb A}_i$, each of
size $<\delta$, such that each cardinal below $\delta$ is collapsed to $\omega$ by 
some ${\mathbb A}_i$, then $\mathbb A$ is isomorphic to $\col(\omega,<\:\delta)$
in such a way that the symmetric reals are exactly those added by the
forcings ${\mathbb A}_i$. In our case the reals added by the forcings ${\mathbb 
B}_0\ast\dots\ast{\mathbb B}_i$ are all in $V$, and eventually all reals of $V$ are 
added. Thus we finally have $\R^V=\dot{\R}_{symm}^{M_\infty}[H]$.

The argument so far is not new. It was first presented by Steel who
used it in \cite{IMMW} to derive several absoluteness results
for $L(\R)$, among them the generic absoluteness of the 
theory of $L(\R)$ with real ---but not ordinal--- parameters. For 
our purposes however this argument is
not sufficient. We have made $\R^V$ the set of reals in the symmetric
collapse of an iterate of $M$, but we must simultaneously make
$\R^{V[G]}$ the set of reals in a different symmetric collapse
of the same iterate. For this reason exactly we left ourselves some
space during the construction, in the form of the embeddings
$j_{2i+1,2i+2}$ and the Woodin cardinals $\delta_{2i+1}$. Let 
$\dot{\mathbb C}_i$ be Woodin's forcing 
$\dot{\mathbb W}^{M,\dot{\mathbb C}_0\ast\dots\ast\dot{\mathbb C}_{i-1}}_{
\delta_{2i},\delta_{2i+1}}$
defined in $M$, and ${\mathbb C}$ their finite support iteration. We will
use those to make the reals of $V[G]$ generic, as we made the
reals of $V$ generic. We must however take care not to spoil the
part of the construction we have completed --- we want to define
$j_{2i+1,2i+2}$ in a way that still allows us to argue that
$\<g(i)\:\mid\:i<\omega\>$ is generic for 
$j_{0,\infty}(\mathbb B)$. For that argument to work we
needed to know that the reals $y_n,\dots,y_{k-1}$ could be chosen in 
$V$, and this followed from
\begin{itemize}
\item[i.] $V^{M_\infty}_{j_{0,\infty}(\delta_i)}$ 
is an element of $V$ for all $i$; and
\item[ii.] $j_{0,\infty}(\delta_i)$ is countable in $V$, for all $i$.
\end{itemize}
Either one of (i),(ii) can easily be maintained using Woodin's first and 
second (see below) genericity iterations. The difficulty is in maintaining both 
conditions simultaneously, and it is here that we must make use of our assumption 
that $P$ is proper. 
\begin{lemma}
\label{lemma}
(Assuming $G$ is $P$--generic/$V$ for some {\em proper} $P$.)
Let $Q=L(V_\delta^Q)$ be uniquely iterable, in $V$, for trees of size 
$\kappa^++1$.
Assume $\delta$ is countable in $V$, let 
$\tau<\eta<\delta$ be ordinals such that
$Q\models\mbox{``$\eta$ is a Woodin cardinal''}$, and consider 
Woodin's forcing ${\mathbb W}={\mathbb W}^Q_{\tau,\eta}$.
{\bf Then} for any real $x\in V[G]$ it is possible to construct an
iteration embedding $j\colon Q\into \tilde{Q}$ 
in $V$ with the property that 
\begin{itemize}
\item $x$ is $j({\mathbb W})$-generic/$\tilde{Q}$;
\item $j(\eta)$ is countable in $V$, indeed $j''(\omega_1^V)\subset
			\omega_1^V$; and
\item $\crit(j)>\tau$.
\end{itemize}
Furthermore for any small forcing ${\mathbb O}\in V_\tau^Q$, if we let
$\dot{\mathbb W}=\dot{\mathbb W}^{Q,\mathbb O}_{\tau,\eta}$ then for any
$o\in V[G]$ which is $\mathbb O$-generic/$Q$ it is possible to construct
an iteration embedding $j\colon Q\into \tilde{Q}$ satisfying the
above except that now $x$ is made $j(\dot{\mathbb W})[o]$-generic
over $\tilde{Q}[o]$.
\end{lemma}

It is worthwhile emphasizing the difference between Woodin's Theorem and Lemma 
\ref{lemma}. In Lemma \ref{lemma}
we allow $x\in V[G]$ (and also $o\in V[G]$ for the second part), 
and still obtain an iteration embedding $j$
in $V$.
Fix $g\colon \omega\into \R^V$ and
$h\colon\omega\into\R^{V[G]}$ so that the pair $g,h$ is
$\col(\omega,\R^V)\cross\col(\omega,\R^{V[G]})$-generic/$V[G]$. 
Granted the Lemma we may 
repeat our construction replacing condition (c) with
\begin{itemize}
\item[c$'$.] $j_{2i+1}\colon M_{2i+1}\into M_{2i+2}$ is an
	embedding 
	to make the real $h(i)$
	generic for the forcing $j_{0,2i+2}(\dot{\mathbb 
C}_i)[h(0),\dots,h(i-1)]$
	over
	the model
	$M_{2i+2}[h(0),\dots,h(i-1)]$. We take $j_{2i+1}$ to be
	the identity if possible. Otherwise we take the embedding given by Lemma 
\ref{lemma}.
\end{itemize}
This modified construction produces $M_\infty$ and $j_{i,\infty}$ satisfying
\begin{enumerate}
\item For all $n<\omega$ 
	$\<g(0),\dots,g(n-1)\>$ is $j_{0,\infty}(\dot{\mathbb B}_0\ast\dots\ast
	\dot{\mathbb B}_{n-1})$-generic/$M_\infty$;
\item For all $n<\omega$ 
      $\<h(0),\dots,h(n-1)\>$ is $j_{0,\infty}(\dot{\mathbb C}_0\ast\dots\ast
	\dot{\mathbb C}_{n-1})$-generic/$M_\infty$; and
\item For $\xi<j_{0,\infty}(\delta)$, $V_\xi^{M_\infty}$
	belongs to $V$ and is countable in $V$.
	\end{enumerate}
Condition (3) and the genericity of $g,h$ allow us as before to argue that
in fact $\< g(i) \mid i<\omega\>$ is $j_{0,\infty}({\mathbb 
B})$-generic/$M_\infty$; and
$\< h(i)\mid i<\omega\>$ is $j_{0,\infty}({\mathbb C})$-generic/$M_\infty$.

As before we can now
convert the forcings ${\mathbb B}$ and 
${\mathbb C}$ into symmetric collapses --- finding $H_1$ and $H_2$ which
are $\col(\omega,<\!j_{0,\infty}(\delta))$-generic/$M_\infty$
so that $\dot{\R}^{M_\infty}_{symm}[H_1]=\so g(i) \mid i<\omega\sr$
and $\dot{\R}^{M_\infty}_{symm}[H_2]=\so h(i) \mid i<\omega\sr$.
Letting $N=M_\infty$ this completes the proof of Theorem \ref{embedding}, at 
least if
$z$ belongs to $M$ --- but if not, before the
beginning of the construction simply iterate
$M$ to make $z$ generic, and then continue to realize $\R^V$ and $\R^{V[G]}$ as 
the reals of a symmetric collapse over $N[z]$.

\vspace{0.1in}
It remains therefore only to prove Lemma \ref{lemma}. We use the following:
\begin{quotethm}[Woodin's second genericity iteration]
Let $Q$ be a $\kappa^++1$-iterable inner model, 
let $\tau<\eta<\kappa^+$ be ordinals
such that $Q\models\mbox{``$\eta$ is a Woodin cardinal''}$, and
let ${\mathbb A}\in V$ be any forcing notion of size $\leq\kappa$. 
Let ${\mathbb W}={\mathbb W}^Q_{\tau,\eta}$ be Woodin's forcing of
the first genericity iteration, defined
in $Q$ from $\tau$ and $\eta$.
{\bf Then} for any $\dot{x}$  which is
a name for a real in $V^{\mathbb A}$, it is possible to construct an iteration
embedding $j\colon Q\into\tilde{Q}$ (in $V$) with the property that
\begin{itemize}
\item For any $F$ which is $\mathbb A$-generic/$V$, the real $\dot{x}[F]$
	is $j({\mathbb W})$-generic/$\tilde{Q}$; 
\item $j(\eta)<(\kappa^+)^V$, indeed $j''(\kappa^+)\subset\kappa^+$; and
\item $\crit(j)>\tau$.
\end{itemize}
Furthermore For any small forcing $\mathbb O\in V^Q_\tau$, if we let 
${\mathbb W}=\dot{\mathbb W}_{\tau,\eta}^{M,{\mathbb O}}$ then for any
$\dot{o}$ which is an $\mathbb A$ name for an $\mathbb O$-generic filter/$Q$, it 
is possible to construct an iteration
embedding $j\colon Q\into \tilde{Q}$ (in $V$) satisfying the
above except that now $\dot{x}[F]$ is made generic over
$\tilde{Q}[\dot{o}[F]]$ (for all $F$ which are $\mathbb A$-generic/$V$).
\end{quotethm}

Using Woodin's second genericity iteration let us prove Lemma 
\ref{lemma}. Let $j^*\colon Q\into Q^*$ be the embedding given
by Woodin's second genericity iteration applied with a name
$\dot{x}$ for the real $x\in V[G]$. Then $j^*\in V$, but $j^*(\eta)$
need not be countable. To overcome this: Fix an elementary
submodel $Y$ of $V_\lambda$ for some sufficiently large $\lambda$
so that $P,\dot{x},Q,j^*,Q^*\in Y$ \footnote{$Q$ is a class model
of course, but it is coded by a real, and we can throw this real into $Y$.}; $Y$ 
belongs to $V$ and is countable in $V$; $G\cap Y$ is $P$-generic/Y; and 
$Y[G\cap Y]\prec V_\lambda[G]$. 
The existence of $Y$ follows from the properness of $P$. In fact
it is enough (by the very definition) to assume that $P$ is reasonable.
Let $\bar{Y}$ be the transitive collapse of $Y$ and $\pi\colon \bar{Y}
\into Y$ the inverse collapse embedding.
Let $\bar{G}={\pi^{-1}}''G$ and 
$\bar{\dot{x}},{j},\tilde{Q}=\pi^{-1}(\dot{x},j^*.Q^*)$.
Let $\bar{x}=\bar{\dot{x}}[\bar{G}]$. Notice that $Q$ is not moved by
$\pi^{-1}$, so we have ${j}\colon Q\into {\tilde{Q}}$.
$\pi$ induces an embedding from $\bar{Y}[\bar{G}]$ onto $Y[G\cap Y]$
which we also call $\pi$. Thus $\pi\colon \bar{Y}[\bar{G}]\into V_\lambda[G]$
is elementary. 

By the elementarity of $\pi$, $\bar{x}$ is 
${j}({\mathbb W}^Q_{\tau,\eta})$-generic/$\tilde{Q}$. Of course
$\bar{x}$ is a real and is not moved by $\pi$,
so $\bar{x}=\pi(\bar{x})=x$.
Thus the embedding ${j}$ makes $x$ generic for Woodin's forcing.
As ${j}\in \bar{Y}$ it is clear that ${j}(\eta)$ is countable in $V$.

The reader can now easily check the remaining requirements
of Lemma \ref{lemma}. Let us here only verify that ${j}$ is an iteration
embedding. This is not obvious --- by elementarity 
$\bar{Y}\models\mbox{``${j}$ is an iteration embedding''}$, but this
does not mean ${j}$ is an iteration embedding in $V$. Let ${\cal T}$
be the iteration tree giving rise to ${j}$. We must show that
the branches ${\cal T}$ chooses are according to the iteration strategy
for $Q$ which we have in $V$. But $Q$ is {\em uniquely} iterable,
so this strategy chooses at every limit stage 
the unique branch with iterable direct limit. Thus it is sufficient to show
that every model $Q_\xi$ on the tree ${\cal T}$ is iterable (in $V$). 
Remember that $\pi$ maps $Q_\xi$ into a model $Q^*_{\xi^*}$ on the tree
${\cal T}^*$ which gives rise to $j^*$. $Q^*_{\xi^*}$ is iterable and
by \cite{it} every model which embeds into an iterable model is
iterable. Thus $Q_\xi$ is iterable and we are done.

The second part of Lemma \ref{lemma} is proved in a similar fashion. Note that
since $Q$ is countable and $\dot{o},{\mathbb O}\in Q$, both
$\dot{o}$ and $\mathbb O$ are automatically in $Y$ and
$\pi(\dot{o})=\dot{o}$. Thus both $\dot{x}[G]$ and $\dot{o}[G]$
are not moved by $\pi\colon \bar{Y}[\bar{G}]\into V_\lambda[G]$.
We take $j^*$ to be the iteration from
the second part of Woodin's second genericity Theorem, and 
immediately by the elementarity
of $\pi$ can conclude that $\dot{x}[G]$ is generic over
$\tilde{Q}[\dot{o}[G]]$. 
\proofend(Lemma \ref{lemma}, Theorem \ref{embedding})

\section{The Anti-Coding Theorem}
\label{sec.anti-coding}

Next let us prove Theorem \ref{anti-coding}. Fix a set $A\subset\ON$ in $V$. We 
must show that $A\in L(\R^V)$ iff
$A\in L(\R^{V[G]})$. Now the implication from left to right
follows immediately from the Embedding Theorem. Assume then that 
$A\in L(\R^{V[G]})$. We must show $A\in L(\R^V)$. As all sets in $L(\R)$
are definable from a real and some ordinals, we may fix a name
$\dot{x}$, ordinals $\vec{\alpha}$ and a formula $\phi$, so
that
$$\check{A}[G]=A=\so \gamma \:\mid\: L(\R^{V[G]})
\models\phi[\vec{\alpha},\dot{x}[G],\gamma]\sr.$$
Without loss of generality we may assume that this is forced by the empty 
condition in $P$.

It is convenient to
replace $A_\kappa$ with the large cardinal assumption $B_\kappa$ stated below. 
It can be seen (using Woodin's second genericity
iteration and some fine structure) that $B_\kappa$ follows from $A_\kappa$. 

\vspace{0.1in}
\noindent
($B_\kappa$) \hspace{0.125in}
For any $K\subset \kappa$ there exists a class model $M$ such that
\begin{itemize}
\item $M=L(V_\delta^M)$, for some $\kappa<\delta<(\kappa^+)^V$;
\item $V_\kappa\subset M$, and $K\in M$;
\item $M\models\mbox{``$\delta$ 
		is the supremum of $\omega$ Woodin cardinals;''}$ and
\item $M$ is uniquely iterable above $\kappa$ for trees of 
	length $\leq(\kappa^+)^V$ (i.e., the good
	player wins the iteration game when the bad player is
	restricted to playing extenders with critical points above
	$\kappa$). 
\end{itemize}
As the forcing $P$ has size $\kappa$ we may take it to be a 
subset of $\kappa$,
and so can fix a model $M$ satisfying the conditions of assumption
$B_\kappa$ with $P,\dot{x}\in M$.
Notice that from $B_\kappa$ it follows that every subset of
$\kappa$ has a sharp, and so $V_\delta^M$ has a sharp.

We now pass to work in a countable elementary submodel $Y\prec V_\lambda$
(for $\lambda$ sufficiently large) which belongs to $V$,
 and contains all relevant objects 
(including $V_\delta^M$ and its sharp).
Let $\bar{Y}$ be the transitive collapse of $Y$, and $\bar{M}$
the image under the collapse map of $M$.\footnote{Again, $M$ is 
a class model. What we mean is that $\bar{M}=L(\bar{V_\delta^M})$
where $\bar{V_\delta^M}$ is the collapse of $V_\delta^M$.}
Let $\pi\colon\bar{Y}\into Y$ be the inverse collapse embedding.
Let $\bar{\dot{x}},\bar{P}$, and $\bar{G}$ be the collapse of
$\dot{x}$, $P$, and $G\cap Y$. Then 
$\bar{\dot{x}}[\bar{G}]=x$, and by properness (reasonability) of
$P$ we may assume that $\bar{G}$ is $\bar{P}$-generic/$\bar{Y}$.
We will attempt to replace the real $\bar{\dot{x}}[\bar{G}]$ in the
definition of $A$ with a real $\bar{\dot{x}}[h]$ for some $h\in V$
which is $\bar{P}$-generic/$\bar{Y}$. The fact that $h\in V$
will then imply that $A\in L(\R^V)$. It is simple to find
$h\in V$ which is $\bar{P}$-generic/$\bar{Y}$ (since $\bar{Y}$ is
countable). The difficulty of course is to do this in such a 
way that $\bar{\dot{x}}[\bar{G}]$ and $\bar{\dot{x}}[h]$ still
define the same set of ordinals.

Let us find $h\in V$ which is $\bar{P}$-generic/$\bar{Y}$ with the
property that for any $\bar{p}\in h$, there exists a condition
$q\leq\pi(\bar{p})$ which is $(Y,P)$-generic.
If $P$ is proper this can be done trivially (perhaps at the price
of modifying $Y$). If $P$ is only known to be reasonable this
is a bit less trivial. Fix in this case some $q_0\in G$ which is
$(Y,P)$-generic. In $V[G]$ there exists an $h$ which is $\bar{P}$-generic/$\bar{
Y}$ such that all conditions in $\pi''(h)$ are compatible with $q_0$
(e.g. take $h=\bar{G}$). By absoluteness then such $h$ exists in $V$,
and it is easy to see that any such $h$ satisfies our requirement above.

Through our choice of $h$ we may, for any condition $\bar{p}\in h$, fix
in some external generic extension of $V$ a filter $G^{\bar{p}}$ 
which is $P$-generic/$V$; contains the condition $\pi(\bar{p})$; and such that 
$\bar{G}^{\bar{p}}={\pi^{-1}}''(G^{\bar{p}}\cap Y)$ is $\bar{P}$-generic over 
$\bar{Y}$ (and hence also over $\bar{M}$ \footnote{We are using here the 
existence of $\bar{V_\delta^M}^\sharp$ inside $\bar{Y}$ to see that all subsets 
of $\bar{P}$ in $\bar{M}$ belong to $\bar{Y}$.}).
By dovetailing together constructions of the sort used in Section 
\ref{sec.embedding} iterate $\bar{M}$ to a model $N$ so that
\begin{itemize}
\item[a.] For each $\bar{p}\in h$ the reals of $V[G^{\bar{p}}]$ 
can be realized as the symmetric collapse over $N[\bar{G}^{\bar{p}}]$; and 
\item[b.] The reals of $V$ can be realized as the symmetric 
collapse over $N[h]$.
\end{itemize}

Let $j\colon \bar{M}\into N$ be the iteration embedding, which we
construct to have critical point above $\bar{\kappa}$, so that
$j(\bar{P})=\bar{P}$, $j(\bar{\dot{x}})=\bar{\dot{x}}$ etc. As in Section 
\ref{sec.embedding} $j$ is a composition of $\omega$ maps, each of which is in 
$V$, and $j$ itself exists only in some external model. Let us denote 
$j(\bar{\delta})$ by $\tilde{\delta}$.\footnote{This is
easily seen to be equal to $\omega_1^V$.} ``The symmetric collapse'' in (a,b) 
above refers to the collapse up to $\tilde{\delta}$.

\begin{claim}
\label{claim.4}
Working in $N[h]$ let $y=\bar{\dot{x}}[h]$ and
consider the forcing $\col(\omega,<\!\tilde{\delta})$.
We claim that for any ordinal $\gamma$ the following are equivalent:
\begin{enumerate}
\item $\gamma\in A$
\item In the forcing $\col(\omega,<\!\tilde{\delta})$ over $N[h]$,
	it is forced that ``$L(\dot{\R}_{symm})\models\phi[\vec{\alpha},
	y,\gamma]$''.
\end{enumerate}
\end{claim}

Otherwise we may fix some ordinal $\gamma$, and a condition 
$\bar{p}\in h$, such that $\gamma\not\in A$ say, and
nonetheless $\bar{p}$ forces in $\bar{P}$ that
$L(\dot{\R}_{symm})\models\phi[\vec{\alpha},
	\bar{\dot{x}},\gamma]$ holds in the
symmetric collapse. (Alternatively $\gamma\in A$
and $\bar{p}$ forces
$L(\dot{\R}_{symm})\not\models\phi[\vec{\alpha},
	y,\gamma]$, but the proof in this case is similar.)
Since $\bar{p}$ is an element of $\bar{G}^{\bar{p}}$ it follows
that over 
$N[\bar{G}^{\bar{p}}]$ it is forced in $\col(\omega,<\!\tilde{\delta})$
that $L(\dot{\R}_{symm})\not\models\phi[\vec{\alpha},
	\bar{\dot{x}}[\bar{G}^{\bar{p}}],\gamma]$.
But now by (a) we may fix $H_{\bar{p}}$
which is $\col(\omega,<\!\tilde{\delta})$-generic/$N[\bar{G}^{\bar{p}}]$
such that $\dot{\R}_{symm}[H_{\bar{p}}]=\R^{V[G^{\bar{p}}]}$.
As furthermore $\bar{\dot{x}}[\bar{G}^{\bar{p}}]=\dot{x}[G^{\bar{p}}]$
it follows that
$L(\R^{V[G^{\bar{p}}]})\not\models\phi[\vec{\alpha},
	\dot{x}[G^{\bar{p}}],\gamma]$.
But this implies $\gamma\not\in A$, a contradiction.
\proofend(Claim \ref{claim.4})

\vspace{0.1in}
By (b) we may fix $H$, a filter which is 
$\col(\omega,<\!\tilde{\delta})$-generic over $N[h]$ so that
$\dot{\R}_{symm}[H]=\R^{V}$. By Claim \ref{claim.4} $\gamma\in A$ iff 
$L(\dot{R}_{symm}[H])\models\phi[\vec{\alpha},y,\gamma]$, i.e.,
$L(\R^V)\models\phi[\vec{\alpha},y,\gamma]$.
Thus, 
$$A=\so \gamma \mid L(\R^V)\models\phi[\vec{\alpha},y,\gamma]\sr\in L(\R^V)$$
completing the proof of the Anti-Coding Theorem.
\proofend(Theorem \ref{anti-coding})

\appendix
\section{Black Boxes}
\label{sec.black}

We include here a proof of Woodin's genericity Theorems. The results
in this Appendix are due to Hugh Woodin 
(circa 1987, to be published in \cite{stf}). 
The reader may easily verify that Woodin's first genericity iteration
is an immediate corollary of the second (taking $\mathbb A$ to be the
trivial forcing for adding nothing and $\kappa=\omega$), 
and so we prove here only
the second. For the rest of this section $\eta$ is assumed 
to be a Woodin cardinal in $Q$.

Consider the algebra ${\cal L}_{\eta}$ of all transfinite formulae
formed by starting with ``$n\in\Ubf{x}$'' (for $n\in\omega$) and
closing under negation and wellordered disjunctions of length $<\eta$.
The forcing ${\mathbb W}_{\tau,\eta}^Q$ is similar to the
Lindenbaum algebra on ${\cal L}_{\eta}$, but rather than
simply setting $\phi \leq \psi \:\iff\: \vdash \mbox{``$\phi\rightarrow
\psi$''}$ Woodin introduces a set of {\em axioms} ${\cal A}\subset 
{\cal L}_{\eta}$ and then defines:
$$\mbox{$\phi\approx\psi$ iff ${\cal A}\vdash \phi
\leftrightarrow \psi$; and $[\phi]\leq [\psi]$ iff ${\cal A}\vdash 
\phi\rightarrow\psi$.}$$
${\mathbb W}_{\tau,\eta}^Q$ is defined to be the forcing notion
consisting of equivalence classes $[\phi]$ for $\phi\in{\cal L}_\eta$,
ordered by $\leq$ as above.

Before writing down the set of axioms ${\cal A}$ note that
with this definition, if $\phi$ is any formula such that $[\phi]=0$
and $x$ any real such that $x\models\phi$, then there must exist an
axiom $a\in \cal A$ such that $x\models\neg a$. Thus any real 
satisfying the axioms cannot satisfy the $0$ condition.

The set ${\cal A}$ is defined as follows: For any $\lambda$ and $\rho$
satisfying
$\tau<\lambda<\rho<\eta$, any $\rho$-strong extender 
$E\in V_\eta^Q$ with $\crit(E)=\lambda$, and any sequence 
$\vec{\phi}=\so\phi_\xi\sr_{
\xi<\lambda}$ of formulae in ${\cal L}_\eta\cap V_\lambda^Q$, 
let $i_E \colon Q\into \ult{Q}{E}$ be the ultrapower embedding
of $Q$, and let $\nu$ be least such that $i_E(\vec{\phi})_\nu$
is not in $V_\rho^Q$. (Notice that $\nu\geq\lambda$, and 
certainly a strict inequality is possible.) The following
formula is taken to be an axiom:
$$a_{\lambda,\rho,E,\vec{\phi}}=\mbox{``$[\bigvee_{\xi<\nu} i_E(\vec{\phi})_\xi]
\: \rightarrow \: [\bigvee_{\xi<\lambda} \phi_\xi]$.''}$$
(This is a formula in ${\cal L}_\eta$, and in fact one which is an
element of $V_{\rho+1}^Q$.) We denote $\nu$ by 
$\nu_{\lambda,\rho,E,\vec{\phi}}$.
It is worthwhile observing that $i_E(\vec{\phi})_\xi=
\phi_\xi$ for $\xi<\lambda$, so that the disjunction $\bigvee_{\xi<\nu} 
i_E(\vec{\phi})_\xi$
is always weaker than (or equal to) the disjunction
$\bigvee_{\xi<\lambda} \phi_\xi$.

Woodin then proves the following Claim
\begin{quoteclaim} In $Q$,
the forcing ${\mathbb W}^Q_{\tau,\eta}$ is $\eta$-c.c.
\end{quoteclaim}
\proof
Assume for contradiction that the Claim fails and fix an anti-chain  
$\so[\psi_\xi]\sr_{\xi<\eta}$ witnessing this. Let $f\colon\eta\into\eta$ be the
function defined by setting $f(\xi)$ to be least $\alpha$ such that
$\psi_\xi\in V_\alpha^Q$. Since $\eta$ is a Woodin cardinal we
can now find $\lambda<\rho$ between $\tau$ and $\eta$ and an extender
$E\in V_\eta^Q$ such that
\begin{enumerate}
\item $\crit(E)=\lambda$;
\item $E$ is $\rho$ strong, and indeed even $\rho$ strong wrt $\so\psi_\xi
	\mid \xi<\eta\sr$; and
\item $\rho>i_E(f)(\lambda)$.
\end{enumerate}
Let $\vec{\phi}=\vec{\psi}\rest\lambda$, and consider the axiom 
$a_{\lambda,\rho,E,\vec{\phi}}$. By condition (3)
$\nu_{\lambda,\rho,E,\vec{\phi}}\geq\lambda+1$ so 
$a_{\lambda,\rho,E,\vec{\phi}}$
clearly proves
$i_E(\vec{\phi})_\lambda\:\rightarrow \bigvee_{\xi<\lambda} \phi_\xi$.
But by condition (2) $i_E(\vec{\phi})_\lambda=\psi_\lambda$
and so ${\cal A}\vdash\mbox{``$\psi_\lambda\rightarrow\bigvee_{\xi<\lambda}
\psi_\xi$''}$. Thus
$[\psi_\lambda]\leq\bigvee_{\xi<\lambda}[\psi_\xi]$ --- a contradiction
since $\so[\psi_\xi]\sr_{\xi<\eta}$ is an anti-chain.
\proofend

\begin{quotelemma}[Woodin] Let $x$ be any real and assume that
$x\models a$ for all $a\in{\cal A}$. Then $x$ generates a ${\mathbb W}_{\tau,
\eta}^{Q}$-generic filter $W_x$, such that $x\in Q[W_x]$.
\end{quotelemma}
\proof 
Define $W_x=\so [\phi]\in {\mathbb W}_{\tau,\eta}^Q \mid x\models \phi\sr$.
This is well defined since $x\models\cal A$.
To see that $W_x$ is a generic filter: Let $\so[\psi_\xi]\sr_{\xi<\beta}$
be a maximal anti-chain in ${\mathbb W}_{\tau,\eta}^Q$ and assume
for contradiction that $[\psi_\xi]\not\in W_x$ for all $\xi<\beta$.
Note that $\beta<\eta$ by the previous Claim, so 
$\varphi=\bigvee_{\xi<\beta}\lnot\psi_\xi$ is a formula in ${\cal L}_\eta$.
$[\varphi]$ is therefore a condition, {\em and} $[\varphi]=0$ since
$\so[\psi_\xi]\sr_{\xi<\beta}$ is a {\em maximal} anti-chain.
But $x\models\varphi$ and this is a contradiction since $x\models{\cal A}$.

Finally to see $x\in Q[W_x]$, note that
$x=\so n \mid [\mbox{``$n\in \Ubf{x}$''}]\in W_x \sr$.
\proofend

\vspace{0.1in}

At last we are in a position to prove Woodin's second genericity Theorem.
Fix a forcing notion $\mathbb A$ of size $\kappa$ and let $\dot{x}$
be a name for a real in $V^{\mathbb A}$. By the previous Lemma, the real
$\dot{x}$ is generic over $Q$ unless it contradicts some of the
axioms in $\cal A$. The reader can easily verify that
if a real $z$ contradicts some axiom $a_{\lambda,\rho,E,\vec{\phi}}$, then
$z$ does 
{\em not} contradict the image axiom $i_E(a_{\lambda,\rho,E,\vec{\phi}})$,
where $i_E\colon Q \into \ult{Q}{E}$ is the ultrapower map.
Thus forming the ultrapower by $E$ ``removes'' the obstruction caused
by the axiom $a_{\lambda,\rho,E,\vec{\phi}}$.
The second genericity iteration is proved by forming an iteration
tree, hitting at every stage the first extender $E$ which defines an
axiom that $\dot{x}$ does not satisfy. A comparison type argument is
then used to show that this iteration terminates. The key to this
comparison type argument is the fact that once an obstructing axiom
has been removed its image will never again become an obstructing axiom. Thus
with each step of the construction we come closer to having no obstructing
axioms at all. This argument requires an iteration {\em tree};
if instead we attempt to use linear iterations then each step may
undo previous steps, and the image of an axiom that was handled previously
may become obstructing again.

Let us begin the construction.
We construct a normal iteration
tree ${\cal T}=\<E_\alpha,\rho_\alpha \mid \alpha<\beta\>$ with
models $Q_\alpha$ and tree structure $T$. The construction is inductive.
At limit $\alpha$ we use our iteration strategy for $Q$ to pick
a cofinal branch of the tree ${\cal T}\rest \alpha$, and set $Q_\alpha$
to be the direct limit of the models along this branch.
At successor stages $\alpha+1$ we must specify $E_\alpha$ and $\rho_\alpha$
(the tree structure is then determined by finding the least $\alpha'$
such that $\rho_{\alpha'}>\crit(E)$ and setting $\alpha' \mathrel{T}\alpha+1$).
We shall use only extenders with critical point above $\tau$.

At successor stages we distinguish between two cases.

\noindent {\bf Case 1}:
If $\dot{x}[F]$ is $j_{0,\alpha}({\mathbb W}_{\tau,\eta}^Q)$-generic/$Q_\alpha$
(for all $\mathbb A$-generic/$V$ filters $F$)
then we let $\beta=\alpha$, $j=j_{0,\alpha}$, and
we are done proving the Theorem.

\noindent {\bf Case 2}: Otherwise, working in $V^{\mathbb A}$ we apply the 
previous
Lemma to $Q_\alpha$ and ${\mathbb W}^{Q_\alpha}_{
\tau,j_{0,\alpha}(\eta)}=j_{0,\alpha}({\mathbb W}_{\tau,\eta}^Q)$,
and conclude that there must be some axiom $a\in j_{0,\alpha}({\cal A})$
such that $\dot{x}\not\models a$. This axiom must have the 
form $a_{\lambda,\rho,E,\vec{\phi}}^{Q_\alpha}$ for some
$\lambda,\rho,E,\vec{\phi}\in Q_\alpha$. Let us pick a condition $q_\alpha\in 
{\mathbb A}$ forcing this, and
forcing value for the unsatisfied axiom $a$, say $q$ forces
$\dot{x}\not\models a_{\lambda_\alpha,\rho_\alpha,E_\alpha,\vec{\phi}^\alpha}$.
Pick $q_\alpha$ so that $\rho_\alpha$ is minimal.
We extend the tree by setting $Q_{\alpha+1}=\ult{Q_{\alpha'}}{E_\alpha}$
for $\alpha'$ least so that $\crit{E_\alpha}=\lambda_\alpha<\rho_{\alpha'}$.

The genericity iteration Theorem will be proved by showing that
the second case in the construction cannot hold for all $\alpha<(\kappa^+)^V$.
This is very similar to the usual proof that comparisons of mice of size
$\kappa$
must terminate before reaching $\kappa^+$.
Assume for contradiction that the construction continues to $(\kappa^+)^V$,
and let $\cal T$ be the tree of length $(\kappa^+)^V$ constructed.
Since $Q$ is assumed to be $(\kappa^+)^V+1$-iterable there exists 
a cofinal branch through the tree. Let $b$ denote this branch. Note that
$b\subset (\kappa^+)^V$ is closed-unbounded.

For every $\alpha\in b$ let $\alpha^+_b$ be the least ordinal such
that $\alpha \mathrel{T} \alpha^+_b+1$. Then $E_{\alpha^+_b}$ has critical
point ($\lambda_{\alpha^+_b}$)
below $\rho_\alpha$, and is applied to $Q_\alpha$ in the tree $\cal T$
to form the ultrapower $Q_{\alpha^+_b+1}$.
Note that $\vec{\phi}^{\alpha^+_b}$ is in $V_{\lambda_{\alpha^+_b}+1}^{
Q_{\alpha^+_b}}$, and since $Q_\alpha$ and $Q_{\alpha^+_b}$
agree on subsets of $\lambda_{\alpha^+_b}$ it follows that
$\vec{\phi}^{\alpha^+_b}\in Q_\alpha$. Let us
denote $\vec{\phi}^{\alpha^+_b}$ by $\vec{\psi}^\alpha$.

Let $S_1\subset b$ be the set of limit points of $b$. For $\alpha\in S_1$
the model $Q_\alpha$ is a direct limit and so
$Q_\alpha=\bigcup_{\beta<\alpha,\beta\in b} j_{\beta,\alpha}''Q_\beta$.
As $\vec{\psi}^\alpha\in Q_\alpha$ there must exist some $h(\alpha)<\alpha$
 such that $\vec{\psi}^\alpha\in j_{h(\alpha),\alpha}'' Q_{h(\alpha)}$.
A standard pressing down argument now produces $\beta<\kappa^+$ and stationary 
$S_2\subset S_1$ so that $h(\alpha)=\beta$ for all $\alpha\in S_2$. Since  
$Q_\beta$ has cardinality $\kappa$, further thinning of $S_2$ produces 
stationary $S_3\subset S_2$ and a fixed $\varphi\in Q_\beta$ such that 
$\vec{\psi}^\alpha=j_{\beta,\alpha}(\vec{\varphi})$
for all $\alpha\in S_3$. Since $\mathbb A$ too has cardinality $\kappa$
we may assume further that for some fixed $q\in \mathbb A$ we have
$q_{\alpha^+_b}=q$ for all $\alpha\in S_3$.

Let $\alpha$ be any element of $S_3$,
and let $\gamma$ be $\alpha^+_b$ (so $\gamma+1\in b$).
Now $q$ forces the real $\dot{x}$ to contradict the axiom
$a^{Q_\gamma}_{\lambda_{\gamma},\rho_{\gamma},E_{\gamma},j_{\beta,\alpha
}(\vec{\varphi)}}$.
This means that necessarily ($q$ forces) $\dot{x}\not\models 
\bigvee_{\xi<\lambda_\gamma} j_{\beta,\alpha}(\vec{\varphi})_\xi$ $(*)$,
{\em and} $\dot{x}\models \bigvee_{\xi<\nu_\gamma}
i_{E_\gamma}^{Q_\gamma}(j_{\beta,\alpha}(\vec{\varphi}))_\xi$.
But $i_{E_\gamma}^{Q_\gamma}(j_{\beta,\alpha}(\vec{\varphi}))=
i_{E_\gamma}^{Q_\alpha}(j_{\beta,\alpha}(\vec{\varphi}))$ 
\footnote{We replaced $Q_\gamma$ with $Q_\alpha$.} since 
$j_{\beta,\alpha}(\vec{\varphi})$ is an element of
$V_{\lambda_\gamma+1}^{Q_\alpha}$. Thus
$\dot{x}\models \bigvee_{\xi<\nu_\gamma}
i_{E_\gamma}^{Q_\alpha}(j_{\beta,\alpha}(\vec{\varphi}))$.
$i_{E_\gamma}^{Q_\alpha}$ is simply $j_{\alpha,\gamma+1}$, so
we can rewrite the above as ($q$ forces) $\dot{x}\models 
\bigvee_{\xi<\nu_\gamma} j_{\beta,\gamma+1}(\vec{\varphi})_\xi$.

Consider now any $\alpha'\in b$ such that $\alpha'>\gamma+1$.
Then $\crit(j_{\gamma+1,\alpha'})\geq\rho_\gamma$ (it is to secure this
fact that we are forced to use iteration trees, and cannot manage with
the simpler linear iterations), and
so for $\xi<\nu_\gamma$, $j_{\beta,\gamma+1}(\vec{\varphi})_\xi$ is
not moved by $j_{\gamma+1,\alpha'}$.
Thus ($q$ forces) $\dot{x}\models \bigvee_{\xi<\nu_\gamma} j_{\beta,\alpha'
}(\vec{\varphi})$.
But then clearly
$\dot{x}\models \bigvee_{\xi<\lambda_{{\alpha'}^+_b}} j_{\beta,\alpha'
}(\vec{\varphi})$, and we now obtain a contradiction to $(*)$
by taking $\alpha'\in S_3$.

This concludes the proof of the first part of the second genericity 
Theorem. We leave the second half to the reader, and indicate
here only how to define the forcing $\dot{\mathbb W}_{\tau,\eta}^{Q,\mathbb O}$
when $\mathbb O$ is a forcing notion in $V_\tau^Q$. Working in
$Q^{\mathbb O}$, again consider the algebra of all formulae
obtained from ``$n\in\Ubf{x}$'' closing under negations and wellordered
disjunctions (in $Q^{\mathbb O}$) of length $<\eta$. Let $\dot{\cal B}$
be the set of axioms (computed in $Q^{\mathbb O}$) $a_{\check{\lambda},
\check{\rho},\dot{E},\dot{\vec{\phi}}}$ as before, with
the restriction that $\dot{E}$ must be an extender (of $Q^{\mathbb O}$)
induced by an extender of $Q$. I.e., there must exist an extender
$F\in Q$ such that the embedding $i_{\dot{E}}^{
Q^{\mathbb O}}\colon Q^{\mathbb O} \into
\ult{Q^{\mathbb O}}{\dot{E}}$ extends the embedding
$i_F^Q\colon Q\into\ult{Q}{F}$.
Set then $$\mbox{$\dot{\phi}\mathrel{\dot{\approx}}\dot{\psi}$ 
iff $\dot{\cal B}\vdash\dot{\phi}\leftrightarrow
\dot{\psi}$; and $\dot{\phi}\mathrel{\dot{\leq}}\dot{\psi}$ iff $\dot{\cal B
}\vdash\dot{\phi}\rightarrow
\dot{\psi}$.}$$
$\dot{\mathbb W}_{\tau,\eta}^{Q, \mathbb O}$ is then defined to be
the set of equivalence classes of $\dot{\approx}$, ordered by $\dot{\leq}$.
The proof of the genericity Theorem proceeds as before. The reader
can verify this, noting that there are many extenders $\dot{E}$
in $Q^{\mathbb O}$ which are induced by extenders in $Q$ --- in fact there
are enough such extenders to witness that $\eta$ is a Woodin cardinal 
(because ${\mathbb O}$ is a ``small'' forcing). This allows carrying
out the argument of the Claim above, and subsequently the
rest of the proof.

\end{document}